%% file: main.tex
\documentclass[11pt, twoside]{article}

\usepackage{amsmath, amssymb, amsthm}
\usepackage{enumitem}
\usepackage[margin=2.71cm]{geometry}
\usepackage{setspace}

\usepackage{dsfont}
\usepackage{relsize}

\usepackage{verbatim}
\usepackage{cases}

\usepackage{graphicx,float,subcaption}
\usepackage{algorithm}
\usepackage{algpseudocode}

\usepackage{hyperref}
\usepackage[svgnames]{xcolor}

\usepackage[numbers,sort,compress]{natbib}

\newtheorem{theorem}{Theorem}

\newtheorem{proposition}[theorem]{Proposition}
\newtheorem{lemma}[theorem]{Lemma}

\theoremstyle{definition}
\newtheorem{assumption}{Assumption}
\newtheorem{definition}{Definition}
\newtheorem{remark}{Remark}

\newcommand{\RR}{\mathbb{R}}

\renewcommand{\Re}{\mathsf{Re}}
\renewcommand{\Im}{\mathsf{Im}}

\newcommand{\EE}{\mathbb{E}}

\newcommand{\norm}[1]{\left\lVert #1 \right\rVert}
\newcommand{\inner}[2]{\left\langle #1, #2 \right\rangle}

\newcommand{\tm}{\mathcal{T}}
\newcommand{\trow}{\mathcal{R}}

\newcommand{\argmin}{\mathop{\mathrm{argmin}}}

\newcommand{\tr}{\mathrm{tr}}

\newcommand{\diag}{\mathrm{diag}}

\newcommand{\opnorm}[1]{\norm{#1}_{\mathrm{op}}}
\newcommand{\hinfnorm}[1]{\norm{#1}_{\mathcal{H}_\infty}}

\newcommand{\bigo}[1]{\mathcal{O}\left(#1\right)}
\newcommand{\bigot}[1]{\widetilde{\mathcal{O}}\left(#1\right)}

\title{\textbf{A least-square method for non-asymptotic identification in linear switching control }}
\author{Haoyuan Sun, Ali Jadbabaie}
\date{}

\begin{document}

\maketitle

\renewcommand{\thefootnote}{\roman{footnote}}
\footnotetext[0]{The authors are affiliated with the Massachusetts Institute of Technology. This work was supported by the ONR grant \#N00014-23-1-2299.}
\footnotetext[0]{Corresponding author: Haoyuan Sun (\texttt{haoyuans [at] mit [dot] edu}).}
\renewcommand{\thefootnote}{\arabic{footnote}}

\begin{abstract}
    The focus of this paper is on linear system identification  in the setting where it is known that the underlying \textit{partially-observed} linear dynamical system lies within a finite collection of known candidate models.
    We first consider the problem of identification from a given trajectory, which in this setting reduces to identifying the index of the true model with high probability.
    We characterize the finite-time sample complexity of this problem by leveraging recent advances in the non-asymptotic analysis of linear least-square methods in the literature.
    In comparison to the earlier results that assume no prior knowledge of the system, our approach takes advantage of the smaller hypothesis class  and leads to the design of a learner with a \textit{dimension-free} sample complexity bound.
    Next, we consider the switching control of linear systems, where there is a candidate controller for each of the candidate models and data is collected through interaction of the system with a collection of potentially destabilizing controllers.
    We develop a dimension-dependent criterion that can detect those destabilizing controllers in finite time.
    By leveraging these results, we propose a \textit{data-driven} switching strategy that identifies the unknown parameters of the underlying system.
    We then provide a non-asymptotic analysis of its performance and discuss its implications on the classical method of estimator-based supervisory control.
\end{abstract}

\section{Introduction}
\input{intro}

\section{Mathematical Preliminaries}
\input{prelim}

\section{Linear System Identification}
\label{sec:linear-identification}
\input{background}

\subsection{System identification from a finite collection}
\label{sec:sample-complexity}
\input{sample-complexity}

\section{Algorithm for identification in switching control}
\label{sec:switching-control}
\input{switching-control}

\section{Conclusion}
\input{conclusion}

\bibliographystyle{plainnat}
\bibliography{reference}

\clearpage

\appendix

\section{Proof of Proposition~\ref{thm:sample-complexity}}
\label{sec:sample-complexity-proof}

\input{sample-complexity-proof}

\section{Proof of Proposition~\ref{thm:detect-unstable}}
\label{sec:detect-unstable-proof}
\input{detect-unstable-proof}

\section{Proof of Theorem~\ref{thm:switching-control}}
\label{sec:switching-control-proof}
\input{switching-control-proof}

\end{document}

%% file: intro.tex
System identification --- the problem of estimating the parameters of an unknown dynamical system from a single trajectory of input/output data --- plays an important role in many problem domains such as control theory, robotics,  and reinforcement learning.
There has been tremendous progress in analyzing the performance of various system identification schemes --- classical results showed asymptotic convergence~\cite{ljung1998system}, whereas recent advances in non-asymptotic theory quantified the sample complexity of learning accurate estimates from data~\cite{oymak2019non,ziemann2023tutorial}.
However, these works all narrowly focus on system identification itself without accounting for the requirements for control applications.
In this work, we consider a problem setting where linear system identification meets switching control so that we develop a \textit{data-driven} approach to simultaneously achieve desirable control and system identification objectives.

In this paper, we consider a collection of discrete-time, partially-observed linear systems $\{(C_i, A_i, B_i)\}_{i=1}^N$ containing the unknown true system parameters $(C_\star, A_\star, B_\star)$. 
These systems are expressed in terms of:
\begin{align*}
    x_{t+1} &= A_i x_t + B_i u_t + w_t, \\
    y_t &= C_i x_t + \eta_t,
\end{align*}
where the dimensions are $x_t \in \RR^{d_x}, u_t \in \RR^{d_u}$ and $y_t \in \RR^{d_y}$.
We assume that the initial state $x_1 \sim \mathcal{N}(0, I_{d_x \times d_x})$, process noise $w_t \sim \mathcal{N}(0, \sigma_w^2 I_{d_x \times d_x})$, and observation noise $\eta_t \sim \mathcal{N}(0, \sigma_\eta^2 I_{d_y \times d_y})$ come from Gaussian distributions.
In many complex systems, e.g. power systems~\cite{meng2016microgrid}, autonomous vehicles~\cite{aguiar2007trajectory}, and public health~\cite{bin2021hysteresis}, it is not practical to design a single controller that achieves satisfactory performance for all candidate models in the collection.
To this end, in a linear switched system, each candidate model has an associated linear controller giving satisfactory performance on this model.
We also note that a mismatched pair of a model and a controller can result in an unstable closed-loop system.

\begin{figure}[!htb]
    \centering
    \includegraphics{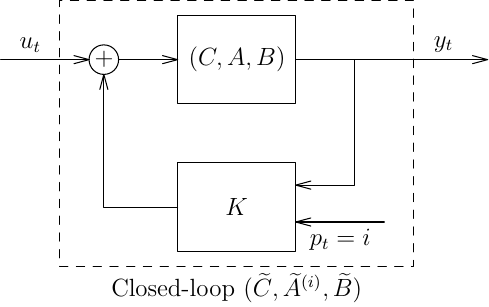}
    \caption{Switching control (with a fixed switching signal)}
    \label{fig:switching-loop}
\end{figure}

Following the convention of \cite{hespanha2001tutorial}, we use the \textit{multi-controller} framework $K(p_t; \check{x}_t, y_t)$, where $\check{x}_t$ is the internal state of the controller, $p_t \in [N]$ is the piece-wise constant switching signal that determines which candidate linear controller is applied, and $y_t$ is the system's output.
And for reasons that will be discussed later, we keep an input signal $u_t$ that is equal to an additive control action on top of the multi-controller.
As illustrated in Figure~\ref{fig:switching-loop}, with an open-loop system $(C, A, B)$ and a fixed switching signal $p_t = j$, the closed-loop system becomes $(\widetilde{C}, \widetilde{A}^{(j)}, \widetilde{B})$, where $\widetilde{A}^{(j)}$ encapsulates both the dynamics $(C, A, B)$ and the controller $K(i; \cdot)$ and $\widetilde{C}, \widetilde{B}$ only depend on $C, B$.
Then, the set of all possible closed-loop dynamics is $\{\{(\widetilde{C}_i, \widetilde{A}^{(j)}_i, \widetilde{B}_i)\}_{i=1}^N\}_{j=1}^N$ and can be pre-computed.
Our goal is then to design a switching strategy that collects the data necessary for identifying the true open-loop parameters $(C_\star, A_\star, B_\star)$ and comes with non-asymptotic performance guarantees.

Although there are several works regarding non-asymptotic identification of partially-observed linear systems~\cite{oymak2019non,simchowitz2019learning,sarkar2021finite,bakshi2023new}, these results all assume that the data are collected from a stable linear system.
As the closed-loop system is potentially unstable in our setting, we cannot directly apply their results.
Furthermore, even under a stabilizing controller, these prior results do not leverage our knowledge that the possible system parameters are contained in a finite set and thus lead to estimation guarantees with sub-optimal dependency on problem dimensions.

Fortunately, there has been extensive work on leveraging the particular properties of switched systems.
Among the existing literature, one popular switching strategy that shares many similarities to our problem setup is the so-called estimator-based supervisory control (see the surveys~\cite{hespanha2001tutorial,liberzon2003switching}).
The estimator-based supervisory control scheme periodically picks the candidate model that most closely matches the observations and applies its associated controller.
While it has been shown that this strategy asymptotically stabilizes the switched system, there are no non-asymptotic guarantees for its performance.
So, we lack a precise characterization of how long this method may take to converge to satisfactory performance.

\paragraph{Contributions.}
In this paper, we focus on the interplay of these two threads of work and derive a novel approach to the study of non-asymptotic system identification in switching control.
To this end, we make the following technical contributions:
\begin{enumerate}[wide, labelindent=3pt, label=\arabic*)]
    \item In Section~\ref{sec:linear-identification}, we present a least-square-based method for linear model identification with prior knowledge that the ground truth is contained in a finite collection of candidate models. Under this setting, we derive a sample complexity bound that is \textit{dimension-free}.
    \item In Section~\ref{sec:instability}, we establish an instability detection criterion by quantitatively bounding the finite-time input-to-output gain of a stable linear system. This allows us to detect any explosive closed-loop dynamics and remove any controllers that are destabilizing the switched system.
    \item Most importantly, in Section~\ref{sec:main-alg}, we present a \textit{data-driven} algorithm for linear system identification problems in switching control. We derive a sample complexity bound on the number of steps for which our strategy finds the correct model with high probability. 
    \item In Section~\ref{sec:implications}, we compare our approach to the classical method of estimator-based supervisory control and discuss the implications of our non-asymptotic guarantees to the problem of switching control.
\end{enumerate}

\subsection{Notations}
For a matrix $M$, we denote $\norm{M}_F$ as its Frobenius norm, $\opnorm{M} = \sigma_{\max}(M)$ as its operator norm (equivalently, its largest singular value), $\rho(M)$ as its spectral radius, and $\tr(M)$ as its trace.
For a stable linear system $(C, A, B)$, we define its $\mathcal{H}$-infinity norm as $\hinfnorm{C, A, B} = \sup_{\norm{s}=1} \sigma_{\max}(C(s I - A)^{-1}B)$, and for simplicity, we use shorthand $\hinfnorm{C, A} = \hinfnorm{C, A, I}$ and $\hinfnorm{A} = \hinfnorm{I, A, I}$.
We also define $P(C, A)$ as the solution to the Lyapunov equation $A^\top P A - P + C^\top C = 0$, and we simply write $P$ when parameters $(C, A)$ are clear from the context.

To simplify our exposition, we sometimes ignore constant factors that do not meaningfully contribute to our conclusions.
We define the big-O notation as $f \in \bigo{g}$ if $\limsup_{x \to \infty} f(x) / g(x) < \infty$ and $f \in \bigot{g}$ if $f \in \bigo{\mathrm{polylog}(\cdot) g(\cdot)}$.
Lastly, we write $f \lesssim g$ if $f \le c \cdot g$ for some universal constant $c$.
Unless otherwise stated, we will explicitly write out any terms dependent on the problem dimensions.
In particular, we consider the Frobenius norm and the trace of a matrix to be dimension-dependent, but the operator norm is not.

\subsection{Literature Review}
Switching control has been studied extensively over the years.
Among the existing literature, there are two popular approaches to designing performant switching policies --- estimator-based supervision that picks the candidate model that most closely resembles the observed process~\cite{hespanha2001tutorial, liberzon2003switching}, and performance-based falsification through some stability certificate~\cite{angeli2002lyapunov,rosa2011stability,stefanovic2008safe,safonov1994unfalsified}.
In this paper, we shall highlight the estimator-based supervision method.
This approach was first formalized in the setting of continuous-time linear switched systems~\cite{morse1996supervisory,morse1997supervisory} and was later extended to nonlinear models~\cite{hespanha1999certainty}, and for discrete-time models~\cite{borrelli1999discrete}.
However, all of the works above only provide asymptotic guarantees for their methods. 

The methods of non-asymptotic linear system identification have seen significant progress with modern tools from statistical learning.
In the case of a fully-observed linear model, \cite{faradonbeh2018finite,simchowitz2018learning} showed that for stable systems, ordinary least square (OLS) achieves estimation error on the order of $\sqrt{T}$, where $T$ is the length of the sample trajectory.
And for unstable linear systems, \cite{sarkar2019near} showed that the OLS estimate may be inconsistent.
Much of the same machinery can be applied to partially-observed stable linear systems, e.g. \cite{oymak2019non,simchowitz2019learning,sarkar2021finite,bakshi2023new}.
Beyond system identification problems, these methods have been applied to problems such as online LQR~\cite{simchowitz2020naive} and latent state learning~\cite{tian2023can}.
A summary of the recent advances in this field can be found in~\cite{ziemann2023tutorial}.

Finally, there are some recent works on applying online learning to switching control.
For example, \cite{li2023online} considers a switched system with fully-observed non-linear models and, inspired by online bandit algorithms, proposes an approach that optimizes for some quadratic cost functions.
We note that the setting of this work differs significantly from ours --- our method applies to partially-observed systems and our ultimate objective is identification so that we do not require access to cost functions.

%% file: prelim.tex
Before we dive into the technical results, we shall briefly introduce some tools from probability and learning theory that are key to our derivations.

We first define a generalization of Gaussian random variables.
Roughly speaking, a sub-Gaussian random variable has tail concentration that is dominated by a Gaussian distribution.
\begin{definition}
    We say that a zero-mean random vector $X \in \RR^d$ is $\sigma^2$-\textit{sub-Gaussian} if for every unit vector $v$ and real value $\lambda$, we have
    \[ \EE \exp(\lambda \inner{v}{X}) \le \exp(\sigma^2\lambda^2 / 2).\]
\end{definition}

In this work, our system identification method is based on ordinary least square (OLS).
Consider a linear model
\[ y_t = \Theta_* z_t + r_t, \, t = 1, \dots, T, \]
where the sequence of random vectors $\{y_t\}_{t=1}^T$ and $\{z_t\}_{t=1}^T$ are adapted to a filtration $\{\mathcal{F}_t\}_{t \ge 0}$ and $r_t$ are residuals/noise that are $\sigma^2$-sub-Gaussian.
Then, OLS seeks to recover the true parameter $\Theta_\star$ through the following estimate:
\[\hat\Theta = \left(\sum_{t=1}^{T} y_t z_t^\top\right) V^{-1}, \text{where } V= \sum_t z_t z_t^\top.\]
We can bound the estimation error $\hat\Theta - \Theta_\star$ with the self-normalized martingale tail bound.

\begin{proposition}[Theorem 1 in \cite{abbasi2011improved}]
\label{thm:self-normalized}
    Consider $\{z_t\}_{t=1}^T$ adapted to a filtration $\{\mathcal{F}_t\}_{t \ge 0}$.
    Let $V = \sum_{t=1}^T z_t z_t^\top$.
    If the scalar-valued random variable $r_t \mid \mathcal{F}_{t-1}$ is $\sigma^2$-sub-Gaussian, then for any $V_0 \succeq 0$,
    \[ \norm{\sum_{t=1}^T z_t r_t}_{(V + V_0)^{-1}}^2 \le 2 \sigma^2 \log\left(\delta^{-1} \frac{\mathrm{det}(V + V_0)^{-1/2}}{\mathrm{det}( V_0)^{1/2}}\right)\]
    with probability $1-\delta$.
\end{proposition}

Additionally, a good OLS estimate requires the covariance matrix $V$ to be non-singular.
Therefore, we want the input data $z_t$'s to be concentrated \textit{away} from 0, which can also be interpreted as persistency of excitation.
We can formalize this concept by the martingale small-ball property.

\begin{definition}[Definition 2.1 in \cite{simchowitz2018learning}]
\label{thm:bmsb-def}
We say that a sequence of $\mathcal{F}_{t}$-adapted random variables $(Z_t)_{t \ge 1}$ satisfies the $(k, v, q)$-\textit{block martingale small-ball (BMSB)} property if for any $t \ge 0$, we have $\frac{1}{k} \sum_{j=1}^k \Pr(Z_{t+j}^2 > v^2 | \mathcal{F}_t) \ge q$ almost surely. 
\end{definition}

This BMSB property implies that the quantity $\sum_{t=1}^T Z_t^2$ scales linearly in $T$ with high probability.
Thus, on average, the random variables $Z_t$ with the BMSB property lie outside of some interval around 0.

\begin{lemma}[Proposition 2.5 in \cite{simchowitz2018learning}]
\label{thm:bmsb-anti-conc}
If the sequence of random variables $(Z_1, \dots, Z_T)$ is $(k, v, q)$-BMSB, then
\[\Pr\left[\sum_{i=1}^T Z_i^2 \le \frac{v^2 q^2}{8} k \lfloor T/k \rfloor \right] \le \exp\left(\lfloor T / k \rfloor q^2 / 8\right).\]
\end{lemma}

Lastly, in this paper, we would encounter the squares of Gaussian random variables.
So, we state a concentration bound on the quadratic form over Guassian random variables.

\begin{proposition}[Corollary 6 in \cite{simchowitz2020naive}]
\label{thm:hanson-wright}
    Consider a symmetric matrix $M \in \mathbb{S}^{d \times d}$ and a random vector $g \sim \mathcal{N}(0, I_{d \times d})$.
    Then, for any $\delta \in (0, 1/e)$,
    \[\Pr\left(|g^\top M g - \tr(M)| > 4 \sigma^2 \tr(M) \log(1/\delta)\right) \le 2 \delta.\]
\end{proposition}
We note this is a simplified version of the Hanson-Wright inequality, and compare to the sharper statement in~\cite{rudelson2013hanson}, we use the fact that $\norm{\cdot}_F$ and $\opnorm{\cdot}$ are bounded by the trace $\tr(\cdot)$.

%% file: background.tex
In this section, we first focus our efforts on system identification.
Specifically, under any constant switching signal $p_t = j$, we identify the index of the true system in the set of possible closed-loop dynamics $\{(\widetilde{C}_i, \widetilde{A}^{(j)}_i, \widetilde{B}_i)\}_{i=1}^N$, which corresponds to the dashed box in Figure~\ref{fig:switching-loop}.
We can then use this index to recover the true parameters of the unknown system.
We stress that the techniques below are applicable to general partially-observed linear systems.
So, in this section, we skip the distinctions between the open-loop vs. closed-loop systems and drop the superscript $\sim$.

\subsection{Problem setup}

We consider a collection of partially-observed linear systems written as:
\begin{align*}
    x_{t+1} &= A_i x_t + B_i u_t + w_t, \\
    y_t &= C_i x_t + \eta_t,
\end{align*}
where the dimensions are $x_t \in \RR^{d_x}, u_t \in \RR^{d_u}$ and $y_t \in \RR^{d_y}$.
We assume that the initial state $x_1 \sim \mathcal{N}(0, I_{d_x \times d_x})$, process noise $w_t \sim \mathcal{N}(0, \sigma_w^2 I_{d_x \times d_x})$, and observation noise $\eta_t \sim \mathcal{N}(0, \sigma_\eta^2 I_{d_y \times d_y})$.
The true system parameters $(C_\star, A_\star, B_\star)$ belong to a collection of $N$ candidate modes $\{(C_i, A_i, B_i)\}_{i=1}^N$.

Additionally, only for this section, we assume that $A_\star$ is stable in the sense that $\rho(A_\star) < 1$.
This assumption is standard in the literature of non-asymptotic linear system identification.
For fully-observed linear systems, \cite{sarkar2019near} gave an example of an unstable linear system that a least-square estimator fails to identify its unknown parameters.
And for the partially-observed setting, as we shall see, it is not possible to bound the residual noise terms when the system is unstable.

There has been a considerable amount of work for the no-prior case where the values of the matrices $(C_\star, A_\star, B_\star)$ can be arbitrary (as long as $A_\star$ is stable), e.g. \cite{oymak2019non, simchowitz2019learning, sarkar2021finite,bakshi2023new}. 
The common approach is to use an exploratory Gaussian noise as the input $u_t \sim \mathcal{N}(0, \sigma_u^2 I_{d_u \times d_u})$.
Then, we use the system's observed response $y_t$ to the sequence of past $h$ Gaussian inputs $z_t := (u_{t-1}, u_{t-2}, \dots, u_{t-h})$ to estimate a Markov parameter that is equal to the system's output controllability matrix with time horizon of length $h$:
\[G_\star := \left[C_\star B_\star, C_\star A_\star B_\star, \dots, C_\star A_\star^{h-1} B_\star\right].\]
Before we apply ordinary least squares (OLS) to estimate the Markov parameter, we first recursively write out the dynamics:
\begin{align*}
    y_t &= C_\star x_t + \eta_t \\
    &= C_\star (A_\star x_{t-1} + B_\star u_{t-1} + w_{t-1}) + \eta_t \\
    &= C_\star A_\star^{h} x_{t-h} \! + \! \underbrace{\sum_{j=1}^h C_\star A_\star^{j-1} B_\star u_{t-j}}_{G_\star z_t} \! + \! \sum_{j=1}^h C_\star A_\star^{j-1} w_{t-j} + \! \eta_t \\
    &:= G_\star z_t + e_t + \eta_t
\end{align*}
Next, we can show that the random vector $e_t = C_\star A_\star^{H} x_{t-H} + \sum_{j=1}^H C_\star A_\star^{j-1} w_{t-j}$
is sub-Gaussian.
For convenience, we write $w_{j:k} = [w_j; w_{j-1}; \dots; w_k]$, $\trow_{k, \ell} = [C_\star A_\star^k, C_\star A_\star^{k+1}, \dots, C_\star A_\star^{k + \ell-1}]$, and $\diag_{m}(B)$ is a block diagonal matrix with $m$ copies of $B$.
We note that $\trow_{k, \ell}$ is a submatrix of the infinite-dimensional Toeplitz operator
\begin{equation}
\label{equ:toeplitz}
\begin{bmatrix} 
    C_\star & C_\star A_\star & C_\star A_\star^{2} & C_\star A_\star^3 & \dots \\
    0 & C_\star & C_\star A_\star & C_\star A_\star^{2} & \dots \\
    0 & 0 & C_\star & C_\star A_\star & \ddots \\
    0 & 0 & 0 & C_\star & \ddots \\
    \vdots & \vdots & \vdots & \ddots & \ddots \\
\end{bmatrix},
\end{equation}
whose operator norm is bounded above by $\hinfnorm{C_\star, A_\star}$ (see Section 4 in \cite{tilli1998singular}).
Therefore, we have $\opnorm{\trow_{k, \ell}} \le \hinfnorm{C_\star, A_\star}$ for any $k, \ell \ge 0$.
Now, we write 
\begin{align*} 
    &C_\star A_\star^h x_{t-H}
    = \trow_{h, t-h-1} w_{t-h-1: 1} + \trow_{h, t-h-1} \diag_{t-h-1}(B_\star) u_{t-h-1: 1} + C A_\star^{t-1} x_1.
\end{align*}
Since $w_{t-h-1: 1}$ is $\sigma_w^2$-sub-Gaussian, the first term is $\opnorm{\trow_{h, t-h-1}}^2 \sigma_w^2 \le \hinfnorm{C_\star, A_\star}^2 \sigma_w^2$-sub-Gaussian.
We then apply this argument to the other terms in $e_t$ and conclude that $e_t$ is $\sigma_e^2$-sub-Guassian for
\[\sigma_e^2 := 2\hinfnorm{C_\star, A_\star} \!\! \sigma_w^2 + \hinfnorm{C_\star, A_\star, B_\star} \!\! \sigma_u^2 + \hinfnorm{C_\star, A_\star}\!.\]
We stress that this argument is only valid when $A_\star$ is stable, otherwise $e_t$ is unbounded.
With this in mind, we express the OLS estimate as
\begin{equation}
\label{equ:ols}
    \hat{G} = \argmin_G \sum_{t=h+1}^{h+\tau} \norm{y_t - G z_t}^2 = \left(\sum_{t=h}^{h+\tau} y_t z_t^\top\right) \Lambda_\tau^{-1}, 
\end{equation}
where $\Lambda_\tau = \sum_t z_t z_t^\top$.

In~\cite{oymak2019non}, it was shown that from a trajectory of length $T$, the OLS estimate $\hat{G}$ satisfies $\opnorm{\hat{G} - G_\star} \le \bigo{\sqrt{h(d_x + d_u) / T}}$.
This bound contains polynomial dependency on the horizon length $h$ and system dimensions $d_x, d_u$ because the size of the Markov parameter $G_\star$ grows with these quantities.
In contrast, this work assumes some prior knowledge that the true parameters comes from a finite set and therefore we proceed to present a least-squares-based approach that yields \textit{dimension-independent} guarantees.

%% file: sample-complexity.tex
We first note that certain collections of candidate models are more difficult to identify than others. 
In particular, more samples would be needed if the collection has a system $(C, A, B)$ whose Markov parameter $G$ is very close to the ground truth $G_\star$.
Thus, we need to quantify how far apart are the models in the given collection in terms of their Markov parameters:
\begin{assumption}
    For all $1 \le i < j \le N$, the Markov parameters in the collection satisfy $\opnorm{G_i - G_j} \ge 2\gamma$.
\end{assumption}

The candidate models would be closer to each other for a smaller value of $\gamma$, which would in turn be harder to distinguish.
Under this assumption, any two candidate models within the collection only need to have different responses to just one input sequence.
Hence, it suffices to come up with estimates that are accurate only respect to these inputs, contrasting to earlier results where the estimation error are uniformly bounded.

As a direct implication of this assumption, for each $1 \le i < j \le N$, there exist unit vectors $u_{ij}, v_{ij}$ so that $|u_{ij}^\top (G_i - G_j) v_{ij}| \ge 2\gamma$.
We call these the \textit{critical directions} of the collection.
It follows that, if an OLS estimate $\hat{G}$ satisfies $|u_{ij}^\top (\hat{G} - G_\star) v_{ij}| < \gamma$ for all $(i, j)$, then the candidate model that is closest to $\hat{G}$ along the critical directions is the ground truth.
This implies that it suffices to find a coarser OLS estimate $\hat{G}$ that is close to the true parameter $G_\star$ in only $\binom{N}{2}$ directions.
We implement this idea as follows:

\begin{algorithm}
\begin{algorithmic}[1]
    \State Input: collection of models $\{G_i\}_{i=1}^N$ and critical directions $\{(u_{ij}, v_{ij}\}_{i < j}$.
    \State Input: Data with $\tau$ samples $\{(y_{H+1}, z_{H+1}), (y_{H+2}, z_{H+2}), \dots, (y_{H+\tau}, z_{H+\tau})\}$.
    \State Compute OLS estimate $\hat{G}$ using \eqref{equ:ols}.
    \State $i \gets 1$
    \ForAll{$j \in \{2, \dots, N\}$}
        \If{$|u_{ij}^\top (G_i - \hat{G}) v_{ij}| \le |u_{ij}^\top (G_j - \hat{G}) v_{ij}|$}
            \State $i \gets j$
            \Comment $j$th model is closer to the estimated $\hat{G}$
        \EndIf
    \EndFor
    \State \Return $i$
\end{algorithmic}
\caption{Linear system model identification with OLS}
\label{alg:ols}
\end{algorithm}

We note that whenever the OLS estimate $\hat{G}$ is accurate along the critical dimensions, the final output given by Algorithm~\ref{alg:ols} must be the index of the correct model, because every other model was shown to not be the closest to $\hat{G}$ along one of the critical directions.
In the following bound, we state the sample complexity of Algorithm~\ref{alg:ols} to identify the correct model with high probability.

\begin{proposition}
\label{thm:sample-complexity}
Given an exploratory input $u_t \sim \mathcal{N}(0, \sigma_u^2 I_{d_u \times d_u})$, there exists a time $\tau \le \frac{\sigma_e^2 + \sigma_\eta^2}{\sigma_u^2 \gamma^2} \log\left(\frac{N^2}{\delta}\right)$, where the inequality holds up to some absolute constants, so that Algorithm~\ref{alg:ols} outputs the correct model with probability $1-\delta$.
\end{proposition}

To interpret this bound, we note that the required number of samples scales quadratically with respect to the ``fineness'' $\gamma$ of the collection.
Also, $\sigma_e^2 + \sigma_\eta^2$ is the magnitude of the residuals and $\sigma_u^2$ is the magnitude of our exploratory inputs, so the quantity $\frac{\sigma_e^2 + \sigma_\eta^2}{\sigma_u^2}$ corresponds to the inverse of the \textit{signal-to-noise ratio} of the system.
In particular, this explains that OLS fails to identify unstable linear systems because in such cases, the signal-to-noise ratio would be arbitrarily low due to unbounded $\sigma_e$.
Finally, the bound only scales logarithmically with respect to $N$, and so the size of the candidate of the collection do not significantly affect the efficacy of our least-square-based method.

\begin{proof}[Proof sketch.]
    Firstly, we can write the estimation error $\hat{G} - G_\star$ purely in terms of the residues $r_t := y_t - G_\star z_t$'s, inputs $z_t$'s and the covariance matrix $\Lambda_t$:
    \[\hat{G} - G_\star = \left(\sum_{t=h+1}^{h+\tau} r_t z_t^\top\right) \Lambda_\tau^{-1}.\]
    Then, for any fixed critical direction $(u, v)$, we apply the self-normalized tail bound (Proposition~\ref{thm:self-normalized}) on the quantity $u^\top (\hat{G} - G_\star) v$ so that, 
    \[\frac{(u^\top (\hat{G} - G_\star) v)^2}{v^\top \Lambda_\tau v + V_0}  \le 2 (\sigma_e^2 + \sigma_v^2) \log \left(\frac{\sqrt{v^\top \Lambda_\tau v + V_0}}{\delta \sqrt{V_0}}\right)\]
    with probability $1-\delta$.
    Next, we want to show persistency of exicitation.
    In particular, if we can pick some value for $V_0$ so that $v^\top \Lambda_\tau v + V_0 \le 2 v^\top \Lambda_\tau v$, then we have a bound on $u^\top (\hat{G} - G_\star) v$.
    To this end, we leverage the block martingale small-ball property (Definition~\ref{thm:bmsb-def}).
    Because $z_t$'s are Guassian, they have a known tail probability and we can show that $v^\top z_t$ is $(1, \sigma_u^2, 3/10)$-BMSB.
    After applying Lemma~\ref{thm:bmsb-anti-conc}, we can show that $v^\top \Lambda_\tau v$ grows linearly in $\tau$ in the sense that there exists some constant $C$ where $\Pr(v^\top \Lambda_\tau v \ge C\tau) \lesssim \exp(-\tau)$.
    These steps lead to a high probability bound on the estimation error along the direction $(u, v)$.

    Finally, we conclude the proof by using union bound over the estimation error along every critical direction $(u_{ij}, v_{ij})_{i < j}$.

    A complete proof can be found in Appendix~\ref{sec:sample-complexity-proof}.
\end{proof}

%% file: switching-control.tex
In this section, we return our focus to the setting of linear switching control.
Recall that in a switched linear system, we assume the unknown underlying linear dynamics is contained in a finite collection of models $\{(C_i, A_i, B_i)\}_{i=1}^N$, and for each model there is an associated linear controller giving satisfactory performance.
As illustrated in Figure~\ref{fig:switching-loop}, the system $(\widetilde{C}_i, \widetilde{A}_i^{(j)}, \widetilde{B}_i)$ represents the closed-loop dynamics when the $j$th controller is applied to the $i$th linear model.
Then, switching control seeks to design a switching strategy that stabilizes the system.
In this work, we additionally want the switching strategy to yield a finite sample guarantee for identifying the unknown system parameters.

\subsection{Summary of estimator-based supervisory control}
One popular approach to switching control is the so-called \textit{estimator-based supervisory control} (see surveys \cite{hespanha2001tutorial, liberzon2003switching}).
At a high-level, this method can be described as follows:
\begin{enumerate}[wide, labelindent=3pt, label=\arabic*)]
    \item We construct a \textit{multi-estimator}, where for at each time $t$ and each model $k$, it takes past outputs $y_t$ and control inputs $u_t$ and makes a prediction $y^{(k)}_{t+1}$ on the next output as if the true underlying system were the $k$th model.
    \item Given the measured outputs $y_t$, we can define the prediction error $\left(e^{(k)}_t := y_t - y^{(k)}_t\right)_{t>0}$.
    \item Let $\hat{i}$ be the index that yields the smallest prediction error according to a time-discounted $\ell_2$-norm and then we apply the $\hat{i}$th controller.
    \item To avoid switching too frequently, we also set a \textit{dwell time} so that we must stick with a switching signal for a prescribed amount of time.
\end{enumerate}

In \cite{morse1996supervisory,morse1997supervisory}, it was shown that the closed-loop switched system resulted from this switching strategy is asymptotically stable in the sense that the system states would remain bounded in response to bounded disturbance.

We note that when the size of the error sequence $e^{(k)}$'s are instead measured by the $\ell_2$-norm without discounting, then this strategy's estimated indexes correspond to the solutions to least-square regression.
This observation motivates us to apply Proposition~\ref{thm:sample-complexity} to derive a switching strategy that has non-asymptotic guarantees.
However, we do not know if the closed-loop dynamics $(\widetilde{C}_i, \widetilde{A}_i^{(j)}, \widetilde{B}_i)$ is actually stable when $i \neq j$.
As we discussed in the previous section, on a trajectory generated by an unstable dynamic, we cannot compute any accurate estimates because the signal-to-noise ratio can be arbitrarily low.
So, before we can determine whether the closed-loop dynamic is stable, running least-square is as good as random guessing.

\subsection{Instability detection}
\label{sec:instability}
The goal of this section is to derive a precise criterion on whether the current closed-loop dynamics is stable, so we can determine if a controller is destabilizing.
As we previously discussed, an unstable partially-observed system would have an arbitrarily low signal-to-noise ratio, which is undesirable for system identification.
So, we exploit the fact that the norms of an unstable system's states would grow without bound.
Under mild assumptions, if the norms of the output are sufficiently large, then we can confidently say that we are facing an unstable system.
Following this intuition, we shall quantify the how explosive are the unstable systems.
\begin{assumption}
    For any unstable closed-loop dynamics $(C_i, A_i^{(j)}, B_i)$, there exists $\varepsilon_a > 0$ so that $\rho(A_i^{(j)}) \ge 1 + \varepsilon_a$.
\end{assumption}

Next, we want to use observability to infer both unstable modes and transient behaviors from past observations.
According to the Hautus (PBH) criterion, a system $(C, A)$ is \textit{observable} if no eigenvector $q$ of $A$ satisfies $Cq = 0$. With this in mind, we make the following assumption.

\begin{assumption}
We say a system $(C, A)$ is \textit{strictly observable} if for every eigenvector $q$ satisfies $\norm{Cq} \ge \varepsilon_c \norm{q}$.
Then, for all $(i, j) \in [N]^2$, $(\widetilde{C}_i, \widetilde{A}_i^{(j)}, \widetilde{B}_i)$ is strictly observable.
\end{assumption}

\begin{remark}
    As a consequence of the Jordan decomposition of $A$, our formulation of strict observability implies the more common definition that $\sigma_{\min}([C; CA; \dots; C A^{d_x-1}]) \ge \varepsilon_c$.
    To see this, we consider an orthonormal basis $\mathcal{B}$ consisting of $A$'s generalized eigenvectors.
    Suppose a basis vector $q \in \mathcal{B}$ is an generalized eigenvector of order $k \le d_x$ and let $q'$ be the unit eigenvector from the same Jordan block.
    Then, $\inner{q'}{A^{k-1}q} = 1$, and from strict observability, we have 
    \begin{align*} 
        \norm{[C; CA; \dots; C A^{d_x-1}] q} \ge \norm{CA^{k-1}q} \ge \left|\inner{q'}{A^{k-1}q}\right| \cdot \norm{Cq'} \ge \varepsilon_c.
    \end{align*}
    Then, for a general vector $v$, we consider its decomposition along this basis and conclude that
    \begin{align*}
        \norm{[C; CA; \dots; C A^{d_x-1}] v}^2
        &= \sum_{q \in \mathcal{B}} \inner{q}{v}^2\norm{[C; CA; \dots; C A^{d_x-1}] q}^2 \\
        &\ge \sum_{q \in \mathcal{B}} \varepsilon_c^2 \inner{q}{v}^2 = \varepsilon_c^2 \norm{v}^2.
    \end{align*}
\end{remark}

These two assumptions together imply that the outputs from an unstable system would be explosive.
In the following result, we employ a threshold \eqref{equ:stable-threshold} corresponding to a high probability bound on the $\ell_2$-norm of the outputs $y_t$ coming from a stable system.
Then, we claim that after a sufficient amount of time, the norms of the outputs exceed this threshold if and only if the dynamics are unstable.

\begin{proposition}
\label{thm:detect-unstable}
    Consider a linear dynamics $(C', A', B')$ belonging to some finite family $\mathcal{S}$. 
    Let the input be $u_t \sim \mathcal{N}(0, \sigma_u^2 I_{d_u \times d_u})$.
    Define
    \begin{equation}
    \label{equ:stable-threshold}
        \begin{aligned}
        \xi(M, \tau, \delta) := \max_{(C, A, B) \in \mathcal{S} \text{ and stable}} \bigg\{2 M 
        + 10 \tau \cdot \left(\sigma_w^2 \tr(P) + \sigma_u^2 \tr(B^\top P B) + \sigma_\eta^2 d_x\right)\log(1/\delta) \bigg\},
        \end{aligned}
    \end{equation}
    where $P$ is the solution to the Lyapunov equation with respect to $(C, A)$.
    Then, each of the following holds with probability $1-\delta$:
    \begin{enumerate}
        \item If the matrix $A'$ is stable and $M \ge x_1^\top P x_1$ over all stable dynamics in $\mathcal{S}$, then  $\sum_{t=1}^\tau \norm{y_t}^2 \le \xi(M, \tau, \delta)$ for any $\tau > 0$.
        \item If the matrix $A'$ is unstable, then $\sum_{t=1}^\tau \norm{y_t}^2 \ge 2 \xi(M, \tau, \delta)$ for some $\tau \lesssim \log((\log(1/\delta)+M)/\delta)$, where the constants have only logarithmic dependency on the dimensions.
    \end{enumerate}
\end{proposition}

We note that the quantity $M$ serves as a bound on the size of the transient. 
And the threshold in \eqref{equ:stable-threshold} in fact is the finite-time input-to-output gain of a stable system.
Therefore, a stable system would not produce outputs that exceed this threshold, whereas the outputs of an unstable system would exceed this threshold after sufficiently large time $\tau$ due to their explosive nature.

\begin{proof}[Proof sketch.]
    The first condition follows from direct computation.
    Specifically, let $[w_{[t-1]}; u_{[t-1]}; \eta_{[t]}]$ be the vector of noises before time $t$.
    Then, the norms of the output $\sum_{s=1}^t \norm{y_s}^2$ can be expressed in terms of a quadratic form
    \[x_1^\top \Sigma_1 x_1 + \bigg[w_{[t-1]}; u_{[t-1]}; \eta_{[t]}\bigg]^\top \mathlarger{\Sigma_2} \bigg[w_{[t-1]}; u_{[t-1]}; \eta_{[t]}\bigg],\]
    where $\Sigma_1$ and $\Sigma_2$ are positive definite matrices whose blocks consist of submatrices from the Toeplitz operator \eqref{equ:toeplitz}.
    We can bound the trace of $\Sigma_1$ and $\Sigma_2$, and then apply the Hanson-Wright inequality (Proposition~\ref{thm:hanson-wright}).
    To interpret the threshold in \eqref{equ:stable-threshold}, we note that the quantity $M$ corresponds to an upper bound on the transient and the second term in \eqref{equ:stable-threshold} corresponds to the steady state response to the process noise and control inputs.
    
    As for the second condition, we note that due to strong observability, the long-term contribution from the process noise $CA^t w_1$ should grow on the order of $(1+\varepsilon_a)^t$.
    Then, we can use this observation to show that $y_t$ in fact satisfies the BMSB property (Definition~\ref{thm:bmsb-def}) with a sufficiently large choice of $\tau$. So, the quantity $\frac{1}{\tau}\sum_{t=1}^\tau \norm{y_t}^2$ scales with $\exp(\tau)$ with high probability.
    Taking logarithm on both sides yields the desired conclusion.

    A complete proof can be found in Appendix~\ref{sec:detect-unstable-proof}.
\end{proof}

As a quick example, we consider the case where $x_1 = 0$ and $M_1 = 0$.
Then, Proposition~\ref{thm:detect-unstable} implies that after some time $\tau \in \bigo{\log(1/\delta) + \log\log(1/\delta)}$, we can tell if the closed-loop dynamics is stable or not.
If we find that we are currently in a stable system, then we can apply the OLS estimation as described in Algorithm~\ref{alg:ols} to determine the true underlying dynamic and apply its corresponding controller for the most desirable performance.

\subsection{Main algorithm and guarantees}
\label{sec:main-alg}
With Propositions~\ref{thm:sample-complexity} and \ref{thm:detect-unstable} in mind, we present our algorithm for switching control that would find the correct model in finite time.

\begin{algorithm}
\begin{algorithmic}[1]
    \State Input: list of dynamics $\mathcal{S} = \{\{(\widetilde{C}_i, \widetilde{A}_i^{(j)}, \widetilde{B}_i)\}_{i = 1}^N\}_{j=1}^N$.
    \State Input: dwell time $\tau_1, \dots, \tau_N$, and $\tau_f$.
    \State Input: upper bound on transients $M_1, \dots, M_\tau$.
    \State We apply exploratory input $u_t \sim \mathcal{N}(0, \sigma_u^2 I_{d_u \times d_u})$.
    \ForAll{$j \in \{1, \dots, N\}$}
        \State Apply $j$th controller for up to $\tau_j$ steps.
        \If{the closed-loop system is stable according to \eqref{equ:stable-threshold} with confidence $1-\frac{\delta}{2N}$.}
            \State Wait for $\bigo{\tau_1 + \dots + \tau_{i-1}}$ steps.
            \State Observe for $\tau_f$ more steps.
            \State Invoke Algorithm \ref{alg:ols} over the collection $\{(C_i, A_i^{(j)}, B_i)\}_{i = 1}^N$ with confidence $1-\delta/2$.
            \State \Return output of Algorithm \ref{alg:ols}
        \EndIf
    \EndFor
\end{algorithmic}
\caption{System identification for switched linear system}
\label{alg:switching-control}
\end{algorithm}

Firstly, on line 4, we use an exploratory input $u_t$ (which we provisioned in Figure~\ref{fig:switching-loop}) to maintain persistency of excitation.
Then, this algorithm works in two stages.
First, on lines 5 -- 7, because the outputs from an unstable dynamics have very little value for learning, we iterate over the list of candidate controllers in some pre-determined order (according to their indexes) and certify their stability with Proposition~\ref{thm:detect-unstable}.
Once we find a controller that leads to a stable closed-loop dynamics, then the results in Section~\ref{sec:linear-identification} are applicable.
Then, on lines 9 -- 11, we roll out a trajectory with the current stable closed-loop system and apply least-square estimation over the set of possible closed-loop dynamics to recover the unknown system parameters.

Before we present the main sample complexity bound for Algorithm~\ref{alg:switching-control}, we first discuss the distinct choices of time $\tau_j$ that we must commit to the $j$th controller.
Recall that the instability detection criterion \eqref{equ:stable-threshold} has two parts: an upper bound on the transient, and the steady-state input-to-output gain from the process and input noises.
The first part depends on the quantity $M$ upper bounding the transient that we must pre-compute.
However, because the controllers may be destabilizing, the internal states of the system are explosive as we apply a greater number of controllers, which leads to larger transients following successive switches.
Therefore, we need to choose $M_j$ that grows with $j$, which in turn requires larger values of $\tau_j$ in order to satisfy the conditions of Proposition~\ref{thm:detect-unstable}.

\begin{theorem}
\label{thm:switching-control}
    We are given appropriate choices on $\tau_f$ according to Proposition \ref{thm:sample-complexity} and $\tau_j = \bigot{j \cdot d_x + \log(1/\delta) d_x}$.
    Then, with probability $1-\delta$, Algorithm~\ref{alg:switching-control} identifies the unknown true system parameters $(C_\star, A_\star, B_\star)$ in 
    \begin{equation}
    \label{equ:switching-total-time}
        \bigot{N^2 d_x + N d_x \log(1/\delta)} + \bigo{\log(N^2/\delta)}
    \end{equation} steps, where the ignored constants have only logarithmic dependency on the system dimensions.
\end{theorem}

There are two components to the sample complexity guarantee in \eqref{equ:switching-total-time}.
The first part is a dimension-dependent term on the time we must take to reject the ``bad'' controllers that are destabilizing.
The second part is a \textit{dimension-independent} sample complexity guarantee on learning the unknown parameters from a stable closed-loop trajectory.
We note that, due to the difficulty in leveraging data produced by an unstable dynamics, our bounds on $M_j$ are very conservative.
This in turn leads to polynomial dependency on both the dimensions and the number of candidate models in the first part of our sample complexity bound.

\begin{proof}[Proof sketch.]
    To analyze Algorithm~\ref{alg:switching-control}, we need to provide upper bound $M_j$ over the transient $\max x_{t_j}^\top P x_{t_j}$, where we denote $t_j$ as the step when the $j$th controller is first applied to the system, and the maximization is taken over all $P$'s that are the solutions to the Lyapunov equations of $(C_i, A^{(j)}_i), i = 1, \dots, N$.
    Finding a suitable $M_1$ is straight-forward, as we recall that $x_1 \sim \mathcal{N}(0, I_{d_x \times d_x})$, and by Proposition~\ref{thm:hanson-wright}, for any positive definite $P$, $x_1^\top P x_1 \le 5 \tr(P) \log(1/\delta)$ with probability $1-\delta$.
    But for subsequent $t_j$'s, the magnitudes of the states $x_{t_j}$'s may grow exponentially quickly as the previous controllers we applied are all destabilizing.
    Nevertheless, using observability, we can carefully bound $x_{t_j}$'s using the past $d_x$ outputs before time $t_j$.
    Note that, if we ensure that $\tau_{j-1} \ge d_x$, we have
    \[ \sum_{t = t_j - d_x}^{t_j-1} \norm{y_t}^2 \le \sum_{t = t_{j-1}}^{t_j-1} \norm{y_t}^2 \le \xi(M_{j-1}, \tau_{j-1}, \delta).\]
    Following a similar direct computation as the first part of Proposition~\ref{thm:detect-unstable}, we can write the quantity 
    \[Y := \sum_{t = t_j - d_x}^{t_j-1} \norm{y_t}^2\]
    in terms of a quadratic form over $x_{t_j - d_x}$ and noise terms.
    After applying Hanson-Wright (Proposition~\ref{thm:hanson-wright}) and the strict observability property, $\norm{x_{t_j - d_x}}^2$ can be upper bounded in terms of $Y$.
    It is worth noting that the constants in this proof are necessarily looser than those of the first part of Proposition~\ref{thm:detect-unstable} because the dynamics are unstable.
    Furthermore, we have $\norm{x_{t_j}}^2 \lesssim \exp(d_x) \norm{x_{t_j - d_x}}^2$.
    It follows that, through an induction argument, we can carefully pick the values of $M_j \lesssim \exp(d_x) M_{j-1}$.
    Finally, through a union bound over all $N$ phases, we can bound the number of steps required to detect those destabilizing controllers.

    Once we find a stabilizing controller, we can leverage our results in Proposition~\ref{thm:sample-complexity}.
    Then, from the current stable closed-loop dynamics, we collect a trajectory whose length is dimension independent and use Algorithm~\ref{alg:ols} to identify the unknown system parameters.

    A complete proof can be found in Appendix~\ref{sec:switching-control-proof}
\end{proof}

\subsection{Implications for estimator-based supervisory control}
\label{sec:implications}

In this section, we discuss the implications of our results for estimator-based supervisory control~\cite{hespanha2001tutorial,liberzon2003switching}.

First, we note that our approach is conceptually quite similar to the estimator-based supervisory control, in that both approaches attempt to determine the model that best describes the unknown system by minimizing the squared-norm of the candidate models' one-step prediction errors against the observed outputs. 
But one major difference is that our approach contains an exploratory and noisy input to ensure persistency of excitation.
This enables us to derive a non-asymptotic sample complexity bound that precisely determines the number of steps we need to take for our least-square estimate to recover the system parameters.
In contrast, the estimator-based supervisory control may not converge to the index of the true model, and thus cannot be used for identification.

Regarding the sample complexity bounds, in Algorithm~\ref{alg:switching-control}, the times $\tau_1, \dots \tau_N$ and $\tau_f$ correspond to the amount of data we must collect to satisfy the conditions of Propositions \ref{thm:sample-complexity} and \ref{thm:detect-unstable}, so that we can learn the model index from the data.
One can view these time values as a precise characterization of  \textit{dwell time}~\cite{morse1996supervisory}.
Under the settings of estimator-based supervisory control, the dwell time is the minimal time interval the switching strategy must commit to a controller before being allowed to switch again.
This dwell time constraint was originally imposed to avoid chattering, but our finite-time analysis endows this quantity with a precise statistical meaning.

Finally, our analysis of instability detection reveals the important role of transient behaviors of the system.
As we previously discussed, the internal states of the system suffer explosive growth from consecutive applications of destabilizing controllers.
Every  time there is a switch to a new controller, the past system states  introduce a transient effect onto the current closed-loop dynamics.
Because the transient would affect the signal-to-noise ratio of the output, it would therefore affect our ability to learn from the data, which results in an increase in the sequence of $\tau_j$'s in Algorithm~\ref{alg:switching-control}.
In particular, our bound in Theorem~\ref{thm:switching-control} indicates that the difficulty in controlling the transient behaviors represent the dominating factor in the sample complexity of Algorithm~\ref{alg:switching-control}.
On the other hand, due to the asymptotic nature of their analysis, existing results on estimator-based supervisory do not take transient terms into account.
One possible solution to this issue would be to use candidate controllers with certain robustness properties so that we are less likely to encounter mismatched pairs of open-loop models and controllers that lead to unstable closed-loop dynamics.

%% file: conclusion.tex
In this paper, we study the problem of non-asymptotic system identification in the context of linear switching control.
We derive a data-driven approach by leveraging ideas from both non-asymptotic system identification and switching control.
In particular, our algorithm works in two stages:
\begin{enumerate}[wide, labelindent=3pt, label=\arabic*)]
    \item We reject any controller that is destabilizing the underlying open-loop dynamics by comparing the observations with our explicit bound on the input-to-output gain of stable systems
    \item Once we certify the stability of closed-loop dynamics, we provide a sharp analysis of system identification that takes into consideration our knowledge of the collection of candidate models.
\end{enumerate}

These ingredients lead to a non-asymptotic guarantee on the sample complexity for learning the unknown system parameters.
From our main results, we reveal new implications on the classical estimator-based supervisory control, particularly regarding to a more precise characterization of the notion of dwell times and the effects of transient behaviors from switching. 

Finally, one future research direction is to derive non-asymptotic guarantees for system identification in \textit{nonlinear} switching control.
Compared to the linear case, the results on the non-asymptotic analysis of nonlinear system identification are significantly more limited.
While it is known that similar guarantees hold for applying ordinary least squares to fully-observed nonlinear systems~\cite{ziemann2022learning}, the partially-observed setting is still an open problem to the best of our knowledge.
Furthermore, translating our Proposition~\ref{thm:detect-unstable} to a nonlinear version seems to be quite difficult because we cannot easily write the outputs purely in terms of the inputs and noise.
So, there many potential works remain in extending the results of this paper to the nonlinear setting.

%% file: sample-complexity-proof.tex
\paragraph{Step 1:}
We start with the definition of the OLS estimate \eqref{equ:ols}:
\begin{align*}
    \hat{G}
    &= \left(\sum_{t=h+1}^{h+\tau} y_t z_t^\top\right) \left(\sum_{t=h+1}^{h+\tau} z_t z_t^\top\right)^{-1} \\
    &= \left(\sum_{t=h+1}^{h+\tau} G_\star z_t z_t^\top\right) \left(\sum_{t=h+1}^{h+\tau} z_t z_t^\top\right)^{-1} + \left(\sum_{t=h+1}^{h+\tau} r_t z_t^\top\right) \left(\sum_{t=h+1}^{h+\tau} z_t z_t^\top\right)^{-1} \\
    &= G_\star + \left(\sum_{t=h+1}^{h+\tau} r_t z_t^\top\right) \Lambda_\tau^{-1},
\end{align*}
where the residual term $r_t$ is equal to $y_t - G_\star z_t = e_t + v_t$, and note that $r_t$ is $(\sigma_e^2 + \sigma_v^2)$-sub-Gaussian.
Therefore, along any critical direction $(u, v)$, we have
\begin{align*}
    (u^\top (\hat{G} - G_\star) v)^2
    &= \left(\sum_{t=h+1}^{h+\tau} u^\top r_t z_t^\top \Lambda_\tau^{-1} v\right)^2. %
\end{align*}

\paragraph{Step 2:}
For this step, we seek to bound the quantity $v^\top \Lambda_\tau v$ for an arbitrary unit vector $v \in \RR^{d_u h}$.
Specifically, for any $\mu > 0$, we want to find an appropriate choice of $\tau$ so that
\begin{equation}
\label{equ:stopping}
    \tau = \min\left\{t \ge 1 : v^\top \left(\sum_{t=h+1}^{\tau+h} z_t z_t^\top\right) v > \frac{2 (\sigma_e^2 + \sigma_v^2)}{\gamma^2} \log\left(\frac{v^\top \left(\sum_{t=h+1}^{h+\tau} z_t z_t^\top\right) v + \mu}{\mu \delta^2 / (9N^4)}\right) \vee \mu\right\}
\end{equation}
with high probability.
Firstly, for the upper bound, we note that $v^\top \Lambda_\tau v = \sum_{t=h+1}^{h+\tau} (v^\top z_t)^2$, which follows a $\chi^2$-distribution. From the concentration bound on $\chi^2$ distributions (Lemma 1, \cite{laurent2000adaptive}), we have
\begin{equation}
\label{equ:cov-upper}
v^\top \Lambda_\tau v \le \sigma_u^2 \tau + 2 \sigma_u^2 \left(\sqrt{\tau \cdot \log(3N^2/\delta)} + \log(3N^2/\delta)\right) \le 2\sigma_u^2 \tau + 3 \sigma_u^2 \log(3N^2/\delta)
\end{equation}
with probability $1-\frac{\delta}{3N^2}$.

As for the lower bound on $v^\top \Lambda_\tau v$, we leverage the ``block martingale small-ball'' (BMSB) property as described in Definition~\ref{thm:bmsb-def}.
Then, since $u^\top z_t$ is $\sigma_u^2$-Gaussian for any unit vector $u$, we have that the input $z_t$'s are $(1, \sigma_u^2 I, 3/10)$-BMSB.
Next, by Proposition 2.5 of \cite{simchowitz2018learning}, we have that
\begin{equation}
\label{equ:cov-lower}
\Pr\left(\sum_{t=h+1}^{h+\tau} (v^\top z_t)^2 \ge \frac{9\sigma_u^2}{800}\tau\right) \le 1 - \exp(9\tau/800).
\end{equation}
Now, we claim that the choice 
\begin{equation}
\label{equ:stopping-bound}
\tau = \max\left(\frac{800}{9}\log(3N^2/\delta), \frac{800\mu}{9\sigma_u^2}, \frac{3200(\sigma_e^2 + \sigma_v^2)}{9\sigma_u^2\gamma^2} \log\left(\frac{16000(\sigma_e^2 + \sigma_v^2) N^4}{\mu \gamma^2 \delta^2}\right) \right)
\end{equation}
satisfies the criteria in \eqref{equ:stopping} with high probability $1-\frac{2\delta}{3N^2}$.

First, we apply \eqref{equ:cov-lower} to conclude that since $\tau \ge \max\left(\frac{800}{9}\log(3N^2/\delta), \frac{800\mu}{9\sigma_u^2}\right)$, we have $v^\top \Lambda_\tau v \ge \mu$ with probability $1 - \frac{\delta}{3N^2}$.
Next, we note that it is not difficult to check that for any positive constant $c$, $x \ge 2c \log (2c)$ implies $x \ge c \log x$.
Therefore, the third term in \eqref{equ:stopping-bound} implies
\begin{align*} 
& \frac{45 \sigma_u^2 N^4}{\delta^2 \mu} \tau \ge \frac{45 \sigma_u^2 N^4}{\delta^2 \mu} \cdot \frac{3200(\sigma_e^2 + \sigma_v^2)}{9\sigma_u^2\gamma^2} \log\left(\frac{45 \sigma_u^2 N^4}{\delta^2 \mu} \cdot \frac{3200(\sigma_e^2 + \sigma_v^2)}{9\sigma_u^2\gamma^2}\right) \\
\implies{}& \tau \ge \frac{1600(\sigma_e^2 + \sigma_v^2)}{9\sigma_u^2\gamma^2} \log\left(\frac{45 \sigma_u^2 N^4}{\delta^2 \mu} \tau\right).
\end{align*}
Recall that from the first term of \eqref{equ:stopping-bound}, $\tau \ge \frac{800}{9}\log(2/\zeta)$.
Therefore, 
\[\tau \ge \frac{1600(\sigma_e^2 + \sigma_v)^2}{9\sigma_u^2\gamma^2} \log\left(\frac{4\sigma_u^2\tau + 6 \sigma_u^2 \log(2/\zeta)}{\delta^2 \mu / (9N^4)} \right).\]
Applying \eqref{equ:cov-upper} and rearranging the terms, we get that
\[\frac{9\sigma_u^2}{800} \tau \ge \frac{2(\sigma_e^2 + \sigma_v)^2}{\gamma^2} \log\left(\frac{2 v^\top \Lambda v}{\delta^2 \mu / (9N^4)} \right) \quad \text{w.p. } 1 - \frac{\delta}{3N}.\]
Using union bound, we conclude that 
\[v^\top \Lambda_\tau v \ge \frac{2(\sigma_e^2 + \sigma_v)^2}{\gamma^2} \log\left(\frac{\mu + v^\top \Lambda v}{\delta^2 \mu / (9N^4)} \right) \vee \mu \quad \text{w.p. } 1 - \frac{2\delta}{3N}.\]

\paragraph{Step 3:}
In this step, we combine the bound we derived from the previous step and self-normalized martingale tail bound.
Note that for any unit vector $u \in \RR^{d_y}$, $u^\top r_t$ is $(\sigma_e^2 + \sigma_v^2)$-sub-Gaussian. 

We now apply Proposition~\ref{thm:self-normalized} on the scalar $(v^\top \Lambda v) z_t^\top \Lambda_\tau^{-1} v = v^\top z_t$ with $V = v^\top \Lambda_\tau v, V_0 = \mu$.
It follows that,
\[\frac{(v^\top \Lambda_\tau v)^2(u^\top (\hat{G} - G_\star) v)^2}{v^\top \Lambda_\tau v + \mu} = \frac{1}{v^\top \Lambda_\tau v + \mu} \left(\sum_{t=h+1}^{h+\tau} (u^\top r_t) (v^\top z_t)\right)^2 \le 2 (\sigma_e^2 + \sigma_v^2) \log \left(\frac{\sqrt{v^\top \Lambda_\tau v + \mu}}{\delta \sqrt{\mu} / (3N^4)}\right)\]
with probability $1-\frac{\delta}{3N^4}$.
Now, when the sample length $\tau$ satisfies \eqref{equ:stopping-bound}, we can apply \eqref{equ:stopping}, so that $v^\top \Lambda_\tau v + \mu \le 2 v^\top \Lambda v$, and
\[(v^\top \Lambda_\tau v)(u^\top (\hat{G} - G_\star) v)^2 \le 4 (\sigma_e^2 + \sigma_v^2) \log \left(\frac{\sqrt{v^\top \Lambda v + \mu}}{\delta \sqrt{\mu} / (3N^4)}\right) = 2 (\sigma_e^2 + \sigma_v^2) \log \left(\frac{v^\top \Lambda v + \mu}{\delta^2 \mu / (9N^4)}\right)\]
with probability $1-\delta/N^2$.
We then apply the condition \eqref{equ:stopping} again to get
\[(u^\top (\hat{G} - G_\star) v)^2 \le \gamma^2 \quad \text{w.p. } 1-\delta/N^2. \]

\paragraph{Step 4:}
With the choice that $\mu = (\sigma_e^2 + \sigma_v^2)/\gamma^2$ and applying union bound over all $\binom{N}{2}$ critical directions $(u_{ij}, v_{ij})_{1 \le i < j \le N}$, we conclude that with
\[\tau = \max\left(\frac{800}{9}\log\left(\frac{3N^2}{\delta}\right), \frac{800(\sigma_v^2+\sigma_e^2)}{9\gamma^2\sigma_u^2}, \frac{3200(\sigma_e^2 + \sigma_v^2)}{9\sigma_u^2\gamma^2} \log\left(\frac{16000 N^4}{\delta^2}\right) \right) \lesssim \frac{\sigma_e^2 + \sigma_v^2}{\sigma_u^2 \gamma^2} \log\left(\frac{N^2}{\delta}\right),\]
the true dynamics is correctly identified with OLS with probability $1-\delta/2$.

%% file: detect-unstable-proof.tex
\paragraph{Part 1:}
For the first part, we a consider linear dynamic $(C, A, B)$ that is stable.
Denote $w_{[t]} = [w_{t}; w_{t-1}; \dots; w_1]$ and similarly for $u_{[t]}$ and $x_{[t]}$.
We have
\[y_{[\tau]} = \tm_{\tau-1} w_{[\tau-1]} + \tm_{\tau-1} \diag_{\tau-1}(B) u_{[\tau-1]} + C[A^{\tau-1}; \dots; A; I] x_1 + \eta_{[t]},\]
where $\diag_{\tau-1}(B)$ is a block-diagonal matrix with $\tau-1$ copies of $B$ and $\tm$ is a Toeplitz matrix
\[ \tm_{\ell} =
\begin{bmatrix} 
    C & CA & CA^{2} &\dots & CA^{\ell-1} \\
    0 & C & C A & \dots & CA^{\ell-2} \\
    0 & 0 & C & \ddots & CA^{\ell-3} \\
    \vdots & \vdots & \ddots & \ddots & \vdots \\
    0 & 0 & \dots & 0 & C \\
    0 & 0 & \dots & 0 & 0
\end{bmatrix}.
\]
It follows that,
\begin{align}
\label{equ:state-quadratic-form}
    & \sum_{t=1}^\tau \norm{y_t}^2
    \le{}& 2 \underbrace{\bigg[w_{[\tau-1]}; u_{[\tau-1]}; \eta_{[\tau]}\bigg]^\top \mathlarger{\Sigma} \bigg[w_{[\tau-1]}; u_{[\tau-1]}; \eta_{[\tau]}\bigg]}_{=: D} 
    + 2 x_1^\top \left(\sum_{t=1}^\tau (A^\top)^{t-1} C^\top C A^{t-1}\right) x_1,
\end{align}
where the covariance matrix $\Sigma$ is
\[ \Sigma = 
\begin{bmatrix}
    \tm_{\tau-1}^\top \tm_{\tau-1} &
    \tm_{\tau-1}^\top \tm_{\tau-1} \diag_{\tau-1}(B) &
    \tm_{\tau-1}^\top I_{\tau d_x \times \tau d_x} \\
    \diag_{\tau-1}(B)^\top \tm_{\tau-1}^\top \tm_{\tau-1} &
    \diag_{\tau-1}(B)^\top \tm_{\tau-1}^\top \tm_{\tau-1} \diag_{\tau-1}(B) &
    \diag_{\tau-1}(B)^\top \tm_{\tau-1}^\top I_{\tau d_x \times \tau d_x} \\
    I_{\tau d_x \times \tau d_x} \tm_{\tau-1} &
    I_{\tau d_x \times \tau d_x} \tm_{\tau-1} \diag_{\tau-1}(B) &
    I_{\tau d_x \times \tau d_x}
\end{bmatrix}\]
Now we bound the quantity $D$.
To apply Hanson-Wright inequality (Proposition~\ref{thm:hanson-wright}), we need to bound trace of $\Sigma$ through direct computation.
Note that 
\begin{align*}
    \tr(\textsf{TM}_{\tau-1}^\top \textsf{TM}_{\tau-1})
    &= \sum_{\ell=0}^{\tau-2} \tr\left(\sum_{k=0}^\ell (A^\top)^k C^\top C A^k\right) \\
    &\le\sum_{\ell=0}^{\tau-2} \tr(P) = \tau \cdot \tr(P).
\end{align*}
Similarly,
\[\tr\left(\diag_{\tau-1}(B)^\top \textsf{TM}_{\tau-1}^\top \textsf{TM}_{\tau-1} \diag_{\tau-1}(B)\right) \le \tau \cdot \tr(B^\top P B).\]
Therefore, by Proposition~\ref{thm:hanson-wright}, we get that for $\delta \in (0, 1/e)$,
\[ \Pr\left(D \ge 5 \tau \left(\sigma_w^2 \tr(P) + \sigma_u^2 \tr(B^\top P B) + \sigma_\eta^2 d_x\right)\log(1/\delta) \right) \le \delta. \]
Finally, since $\sum_{t=1}^\tau (A^\top)^{t-1} C^\top C A^{t-1} \preceq P$, we conclude that

\begin{equation*}
\label{equ:stable-upper}
    \sum_{t=1}^\tau \norm{y_t}^2 \le 2x_1^\top P x_1 + 10 \tau \cdot \left(\sigma_w^2 \tr(P) + \sigma_u^2 \tr(B^\top P B) + \sigma_\eta^2 d_x\right)\log(1/\delta) \quad \text{w.p. } 1-\delta.
\end{equation*}

\paragraph{Part 2:} 

For a strictly unstable $A$, we claim that strict observability ensures $\opnorm{CA^t} \ge \frac{1}{\sqrt{2}} \varepsilon_c (1+\varepsilon_a)^t$ for any $t > 0$.
To see this, let $q$ be an eigenvector of $A$ with eigenvalue $\lambda$.
Then, we have $\norm{C A^t q} \ge |\lambda|^t \varepsilon_c$.
WLOG, let us assume $\norm{\Re(C A^t q)} \ge \norm{\Im(C A^t q)}$, which means $\norm{C A^t \Re(q)}^2 \ge \frac{1}{2} |\lambda|^{2t} \varepsilon_c^2$.
Since $\rho(A) \ge 1+\varepsilon_a$, we have $\opnorm{CA^t} \ge \frac{1}{\sqrt{2}} \varepsilon_c (1+\varepsilon_a)^t$.

Thus, for a fixed $\tau$, there exists unit vectors $(u, v)$ so that $|u^\top CA^\tau v| \ge \varepsilon_c (1+\varepsilon_a)^\tau$.
Then, for any $t > \tau$, we have
\begin{align*}
    \EE[(u^\top y_t)^2] 
    &= \EE[u^\top y_t y_t^\top u] \\
    &\ge \EE[u^\top (CA^\tau w_{t-\tau}) (CA^\tau w_{t-\tau})^\top u] \\
    &\ge \EE[u^\top (CA^\tau v v^\top w_{t-\tau}) (CA^\tau v v^\top w_{t-\tau})^\top u] \\
    &\ge \varepsilon_c (1+\varepsilon_a)^\tau \EE[(v^\top w_{t-\tau})^2] = \varepsilon_c (1+\varepsilon_a)^\tau \sigma_w^2.
\end{align*}
Then, through the Paley-Zygmund lower bound (see (3.12) in \cite{simchowitz2018learning}), we find that the sequence $(u^\top y_t)_{t> \tau}$ is $(1, \sqrt{\varepsilon_c (1+\varepsilon_a)^\tau} \sigma_w, 3/10)$-BMSB.
Hence, using Proposition 2.5 in \cite{simchowitz2018learning}, we get that
\begin{equation*}
\sum_{t=1}^{2\tau} \norm{y_t}^2 \ge \sum_{t=\tau+1}^{2\tau} (u^\top y_t)^2 \ge \frac{9\varepsilon_c (1+\varepsilon_a)^\tau \sigma_w^2}{800} \tau \quad \text{w.p. } 1 - e^{-9\tau/800}.
\end{equation*}
From here, it is not difficult to check that with a choice of 
\begin{equation}
\label{equ:unstable-thres}
\tau =  \max\left\{\frac{1600}{9} \log\left(\frac{1}{\delta}\right), \frac{1}{\log(1+\varepsilon_a)} \log \left(\frac{6400 \varepsilon_c}{9\sigma_w^2\delta} \left(M + 5 (\sigma_w^2 \tr(P) + \sigma_u^2 \tr(B^\top P B) + \sigma_\eta^2 d_x)\log\left(\delta^{-1}\right)\right)\right)\right\},
\end{equation}
whenever $A$ is unstable, we have $\sum_{t=1}^{\tau} \norm{y_t}^2 \ge 2 \xi(M, \tau \delta)$ with probability $1-\delta$.

%% file: switching-control-proof.tex
\paragraph{Part 1.}
For convenience, we denote the set all possible closed-loop dynamics $\{\{\widetilde{C}_i, \widetilde{A}^{(j)}_i, \widetilde{B}_i\}_{i=1}^N\}_{j=1}^N$ as $\mathcal{S}$, $\delta' = \delta/(4N)$, and $\bar\xi_j = \xi(M_{j}, \tau_{j}, \delta')$.
Furthermore, we define some constants expressed in terms of the problem parameters:
\begin{align*}
    m_a &:= \max_{(C, A, B) \in \mathcal{S}} \max\{1, \opnorm{A}\},\\
    m_s &:= \max_{(C, A, B) \in \mathcal{S}} \max\{1, \opnorm{B}, \opnorm{C}\},\\
    m_p &:= \max_{(C, A, B) \in \mathcal{S} \text{ and stable}} \opnorm{P}, \text{where } A^\top P A - P + C^\top C = 0,\\
    m_t &:= \max_{(C, A, B) \in \mathcal{S} \text{ and stable}} \tr(P),\\
    \sigma_m &:= \max\{\sigma_w, \sigma_u, \sigma_\eta\},\\
    c_e &:= \max\{1, 1/\log(1+\varepsilon_a)\},\\
    c_r &:= m_p (22 \varepsilon_c^{-2} + 1)\sigma_m^2 c_e,\\
    c_p & := 2 \cdot \max\{1, m_p \varepsilon_c^{-2}\}, \\
    c_s &:= \max_{(C, A, B) \in \mathcal{S} \text{ and stable}} \left\{5 \left(\sigma_w^2 \tr(P) + \sigma_u^2 \tr(B^\top P B) + \sigma_\eta^2 d_x\right)\log(1/\delta') \right\}.
\end{align*}

Through induction, we show that under appropriate choice of 
\[\tau_j = (j-1) \left(\frac{2}{\log(m_a)} d_x + \log(c_p)\right) + \tau_1,\]
where $\tau_1$ is defined in \eqref{equ:tau_1}, and some recursively defined $M_j$ that we shall reveal later,
we have
\begin{enumerate}[wide, labelindent=2pt,label=\textbf{\arabic*})]
    \item $\frac{9\varepsilon_c (1+\varepsilon_a)^{\tau_{j}} \sigma_w^2}{800} \tau_{j}  \ge \bar\xi_j + c_r \log(2/\delta') d_x^2 m_a^{4d_x}$;
    \item for $P$ that is the solution of Lyapunov equation of some stable $(C, A)$ in $\mathcal{S}$, we have $x_{t_j}^\top P x_{t_j} \le M_j$ with probability $1-\delta'$.
\end{enumerate}

The base case where $i = 1$ is simple.
Recall that $x_1 \sim \mathcal{N}(0, I_{d_x t \times d_x})$, then by Proposition~\ref{thm:hanson-wright}, $x_1^\top P x_1 \le 5 m_t \log(1/\delta')$ with probability $1-\delta'$.
Then, from \eqref{equ:unstable-thres}, it suffices to take $M_1 = 5 m_p \log(1/\delta')$ and $\tau_1 = \bigo{d_x + \log(1/\delta')}$, where
\begin{equation}
\label{equ:tau_1}
\tau_1 =  \max\left\{\frac{1600}{9} \log\left(\frac{1}{\delta'}\right), \frac{1}{\log(1+\varepsilon_a)} \log \left(\frac{6400 \varepsilon_c}{9 \sigma_w^2\delta'} \left(M_1 + c_r \log(2/\delta') d_x^2 m_a^{4d_x} + c_s \log(1/\delta')\right)\right)\right\}
\end{equation}
so that we can tell if the first controller is destabilizing or not with probability $1-\delta/(2N)$.

Now assuming that we have suitable choice of $\tau_{j-1}$, we show that we can tell if the $j$th controller is destabilizing or not by taking $\tau_j = \bigo{d_x} + \tau_{j-1}$.
For simplicity, we denote the $(j-1)$st system as $(C, A, B)$.
Note that, by following similar steps as we did in \eqref{equ:state-quadratic-form}, we have:
\begin{align*}
    & x_{t_j-d_x}^\top \left(\sum_{t=1}^\tau (A^\top)^{t-1} C^\top C A^{t-1}\right) x_{t_j-d_x}  \\
    \le{}& 2 \sum_{t=t_j-dx}^{t_j-1} \norm{y_t}^2 + 2 \underbrace{\bigg[w_{t_j-1:t_j-d_x}; u_{t_j-1:t_j-d_x}; \eta_{t_j:t_j-d_x}\bigg]^\top \mathlarger{\Sigma} \bigg[w_{t_j-1:t_j-d_x}; u_{t_j-1:t_j-d_x}; \eta_{t_j:t_j-d_x}\bigg]}_{=: D} 
\end{align*}
where the covariance matrix $\Sigma$ is
\[ \Sigma = 
\begin{bmatrix}
    \tm_{d_x-1}^\top \tm_{d_x-1} &
    \tm_{d_x-1}^\top \tm_{d_x-1} \diag_{d_x-1}(B) &
    \tm_{d_x-1}^\top I_{d_x^2 \times d_x^2} \\
    \diag_{d_x-1}(B)^\top \tm_{d_x-1}^\top \tm_{d_x-1} &
    \diag_{d_x-1}(B)^\top \tm_{d_x-1}^\top \tm_{d_x-1} \diag_{d_x-1}(B) &
    \diag_{d_x-1}(B)^\top \tm_{d_x-1}^\top I_{d_x^2 \times \tau d_x^2} \\
    I_{\tau d_x \times \tau d_x} \tm_{\tau-1} &
    I_{\tau d_x \times \tau d_x} \tm_{\tau-1} \diag_{\tau-1}(B) &
    I_{\tau d_x \times \tau d_x}
\end{bmatrix}\]
Then, we bound $D$ with Hanson-Wright inequality (Proposition~\ref{thm:hanson-wright}).
Note that 
\begin{align*}
    \tr(\tm_{d_x-1}^\top \tm_{d_x-1})
    &= \sum_{\ell=0}^{d_x-2} \tr\left(\sum_{k=0}^\ell (A^\top)^k C^\top C A^k\right) \\
    &\le \sum_{\ell=0}^{d_x-2} \opnorm{C^\top C} \sum_{k=0}^\ell d_x \opnorm{A^k (A^\top)^k } \\
    &\le \sum_{\ell=0}^{d_x-2} d_x \opnorm{C}^2 \frac{\opnorm{A}^{2\ell+2}}{\log \opnorm{A}} \\
    &\le d_x^2 \cdot m_s^2 \frac{m_a^{2d_x}}{\log(1+\varepsilon_a)}.
\end{align*}
Similarly,
\[\tr\left(\diag_{d_x}(B)^\top \tm_{d_x-1}^\top \tm_{d_x-1} \diag_{\tau-1}(B)\right) \le d_x^2 \cdot m_s^4 \frac{m_a^{2d_x}}{\log(1+\varepsilon_a)}.\]
Therefore, by Proposition~\ref{thm:hanson-wright}, we get that for $\delta \in (0, 1/e)$,
\[ \Pr\left(D \ge 5 d_x^2 \left((\sigma_w^2 + \sigma_u^2) m_a^{2d_x} \log(1+\varepsilon_a)^{-1} + \sigma_\eta^2\right)\log(2/\delta') \right) \le \delta'/2. \]
Because of
\[ \sum_{t = t_j - d_x}^{t_j-1} \norm{y_t}^2 \le \sum_{t = t_{j-1}}^{t_j-1} \norm{y_t}^2 \le \xi(M_{j-1}, \tau_{j-1}, \delta') = \bar\xi_{j-1},\]
and strict observability,
we conclude that
\[\norm{x_{t_j-d_x}}^2 \le \varepsilon_c^{-2} x_{t_j-d_x}^\top \left(\sum_{t=1}^\tau (A^\top)^{t-1} C^\top C A^{t-1}\right) x_{t_j-d_x} \le 2\varepsilon_c^{-2} \bar\xi_{j-1} + 22 \varepsilon_c^{-2} \sigma_m^2 c_e d_x^2 m_a^{2d_x} \log(2/\delta')
\]
with probability $1-\delta'/2$. Next, we can write
\[x_{t_j} = A^{d_x} x_{t_j - d_x} + \sum_{k=1}^{d_x} A^{k-1} w_{t_j - k} + \sum_{k=1}^{d_x} A^{k-1} B u_{t_j - k}.\]
Therefore, by applying Proposition~\ref{thm:hanson-wright}, we have
\begin{align*}
    \norm{x_{t_j}}^2 
    &\le 3 \norm{A^{d_x} x_{t_j - d_x}}^2 + 3\norm{\sum_{k=1}^{d_x} A^{k-1} w_{t_j - k}}^2 + 3\norm{\sum_{k=1}^{d_x} A^{k-1} B u_{t_j - k}}^2 \\
    &\le 3 \opnorm{A}^{2d_x} \norm{x_{t_j - d_x}}^2 + 15 \log(2/\delta') \left(\sigma_w^2 \tr\left(\sum_{k=1}^{d_x} (A^\top)^{k-1} A^{k-1}\right) + \sigma_u^2 \tr\left(\sum_{k=1}^{d_x} B^\top (A^\top)^{k-1} A^{k-1} B\right) \right) \\
    &\le 3 m_a^{2d_x} \norm{x_{t_i - d_x}}^2 + \bigo{\log(2/\delta') d_x^2 \opnorm{A}^{4d_x}} + 15 (\sigma_w^2 + \sigma_u^2) \log(2/\delta') d_x \opnorm{A}^{2d_x} \\
    &\le m_a^{2d_x} \varepsilon_c^{-2} \bar\xi_{i-1} + (22 \varepsilon_c^{-2} + 1)\sigma_m^2 c_e d_x^2 m_a^{4d_x} \log(2/\delta')
\end{align*}
with probability $1-\delta'$.
And it follows that for any $P$ that is the solution of Lyapunov equation of some stable $(C, A)$ in $\mathcal{S}$, we have  \[ x_{t_j}^\top P x_{t_j} \le \opnorm{P} \norm{x_{t_i}}^2 \le m_p \varepsilon_c^{-2} m_a^{2d_x} \bar\xi_{i-1} + c_r \log(2/\delta') d_x^2 m_a^{4d_x}\]
with high probability.

With the choice of $M_j = m_p \varepsilon_c^{-2} m_a^{2d_x} \bar\xi_{i-1} + c_r \log(2/\delta') d_x^2 m_a^{4d_x}$, we note that
\begin{align*}
    \bar\xi_j 
    = M_j + \tau_j \cdot c_s 
    & \le m_p \varepsilon_c^{-2} m_a^{2d_x} \bar\xi_{j-1} + c_r \log(2/\delta') d_x^2 m_a^{4d_x} + \tau_j \cdot c_s \\
    &= m_p \varepsilon_c^{-2} m_a^{2d_x} M_{j-1} + c_r \log(2/\delta') d_x^2 m_a^{4d_x} + (m_p \varepsilon_c^{-2} m_a^{2d_x} \tau_{j-1} + \tau_{j}) c_s \\
    &\le  m_a^{2d_x} \left(m_p \varepsilon_c^{-2} M_{j-1} + c_r \log(2/\delta') d_x^2 m_a^{2d_x} + (m_p \varepsilon_c^{-2} + 1) \tau_j c_s \right) \\
    &\le  c_p m_a^{2d_x} \left(M_{j-1} + c_r \log(2/\delta') d_x^2 m_a^{2d_x} + \tau_j c_s \right)
\end{align*}
Note that by the inductive hypothesis, 
\[\frac{9\varepsilon_c (1+\varepsilon_a)^{\tau_{j-1}} \sigma_w^2}{800} \tau_{j-1}  \ge M_{j-1} + c_r \log(2/\delta') d_x^2 m_a^{2d_x} + \tau_{j-1} c_s.\]
Therefore, it suffices to take $\tau_j = \tau_{j-1} + \frac{2}{\log(m_a)} d_x + \log(c_p)$ so that
\[\frac{9\varepsilon_c (1+\varepsilon_a)^{\tau_{i}} \sigma_w^2}{800} \tau_{i}  \ge \bar\xi_i + c_r \log(2/\delta') d_x^2 m_a^{2d_x}.\]
By union bound and Proposition~\ref{thm:detect-unstable}, we can correctly determine whether the $j$th controller is stabilizing with probability $1-\delta/(2N)$.
Then, we can correctly determine whether any of the $N$ controllers is destabilizing with probability $1-\delta/2$.

\paragraph{Part 2.}
The algorithm becomes straightforward after we find a stabilizing controller.
Note that we spend up to $\sum_{j=1}^N \tau_j \in \bigot{N^2 d_x + N d_x \log(1/\delta')}$ steps, where we ignore constants that are logarithmic in problem dimensions, on those destabilizing controllers before finding a suitable stabilizing controller.
To apply our least-square identification result, we first need to maintain the current stabilizing controller until the states are bounded, e.g. $\norm{x_t}^2 \le 1$, and this takes another $\bigot{N^2 d_x + N d_x \log(1/\delta')}$ steps.

Finally, we directly apply Algorithm~\ref{alg:ols} and Proposition~\ref{thm:sample-complexity}.
And after $\tau_f \lesssim \frac{\sigma_e^2 + \sigma_\eta^2}{\sigma_u^2 \gamma^2} \log\left(\frac{4N^2}{\delta}\right)$ steps, we find the index of the true system with probability $1-\delta$.
Overall, Algorithm~\ref{alg:switching-control} recover the true ssytem parameters with probability $1-\delta$.